\documentclass[leqno,11pt]{amsart}
  \usepackage{amsmath}
  \usepackage{amssymb,amsxtra,amsthm,amscd}
  \usepackage{blaom}
  \newcommand{\nabladot}{\nabla_{\text{\Large$\cdot$}\,}}
  \newcommand{\gl}{\mathfrak{gl}}
  \DeclareMathOperator{\operator}{oper}
  \newcommand{\levi}{\nabla^{\mathrm{L}\text{-}\mathrm{C}}}
\begin{document}
\bibliographystyle{plain}
\title[Geometric structures]%
{Geometric structures\\ as deformed infinitesimal symmetries}
\author{Anthony D.~Blaom}%
\date{\today}%
\address{Department of Mathematics, University of Auckland, Private
  Bag 92019, Auckland, New Zealand.}%
\thispagestyle{empty}
\begin{abstract}
  A general model for geometric structures on differentiable
  manifolds is obtained by deforming infinitesimal symmetries.
  Specifically, this model consists of a Lie algebroid, equipped with
  an affine connection compatible with the Lie algebroid structure.
  The curvature of this connection vanishes precisely when the
  structure is locally symmetric.
  
  This model generalizes Cartan geometries, a substantial class, to
  the intransitive case.  Simple examples are surveyed and
  corresponding local obstructions to symmetry are identified.  These
  examples include foliations, Riemannian structures, infinitesimal
  $G$-structures, symplectic and Poisson structures.
\end{abstract}
\maketitle
\section{Introduction}
According to \'Elie Cartan, a geometric structure is a symmetry
deformed by curvature.  Here we describe a model for geometric
structures promoting this vision, formulated in the language of Lie
algebroids.

If by {\df symmetry} we mean a smooth action $G_0\times M\rightarrow
M$ of a Lie group $G_0$ on a smooth manifold $M$, then every symmetry
has an infinitesimal counterpart: the corresponding action $\mathfrak
g_0\times M\rightarrow TM$ of the Lie algebra $\mathfrak g_0$ of
$G_0$.  Such infinitesimal actions are generalized by vector bundle
morphisms $\mathfrak g\rightarrow TM$, where $\mathfrak g$ is a
possibly non-trivial vector bundle over $M$ known as a {\df Lie
  algebroid}. We call a Lie algebroid equipped with an affine
connection $\nabla$ a {\df Cartan algebroid} whenever $\nabla$ is
compatible with the Lie algebroid structure in an appropriate sense.
As it turns out, one may then view the curvature of such a connection as
the local obstruction to symmetry.  

The significance of Cartan algebroids is that they are a natural
model for many geometric structures.  In a sense, they are
infinitesimal versions of Cartan's {\df espace generalis\'e}, also
known as {\df Cartan geometries}.  Cartan geometries include first
order structures, such as Riemannian and almost Hermitian structures,
and higher order structures, such as projective, conformal and CR
structures.  For some details and further examples, see
\cite{Sharpe_97} and \cite{Slovak_97}.

The present generalization adds {\em intransitive} structures
(deformations of intransitive symmetries) to the list. These include
foliations and Poisson manifolds, equipped with suitable connections.
Even in the transitive case however, the Cartan algebroid point of
view is somewhat novel.

For instance, `curvature' assumes a new meaning.  Classical
notions of curvature measure deviation from a {\em particular}
symmetric model.  In a Cartan geometry, this model is a prescribed
homogeneous space $G/H$. In the context of $G$-structures (see, e.g.,
\cite{Kobayashi_72}), the implicit model is usually $\mathbb R^n$.

By contrast, the present theory {\em has no models:} all symmetric
structures are created equal and the curvature of $\nabla$ merely
measures deviation from {\em some} symmetric structure. For example,
Euclidean space, hyperbolic space and spheres are all regarded by us
as `flat' Riemannian structures.  

Beyond the simplest of examples, however, the procedure by which a
Cartan algebroid is associated with a given geometric structure is not
trivial. This procedure amounts to a model-free version of Cartan's
method of equivalence. For the $G$-structure implementation of this
method, see, e.g., \cite{Gardner_89}. A Lie algebroid version of
Cartan's method is developed in detail in our paper \cite{Blaom_A},
which also refines or improves several ideas introduced in the
present work.


\subsection*{Paper outline}
In Sect.\,\ref{section2} we recall how Lie algebroids may be viewed as
generalized Lie algebra actions. We then formulate what it means for
an affine connection on a Lie algebroid to be compatible with the Lie
algebroid structure (Definition \ref{geom}). Such connections, which
we call {\df Cartan connections}, are related to the classical
connections of the same name (see Sect.\,\ref{section7}). Theorem A
(Sect.\,\ref{section2}) which characterizes the locally symmetric
Cartan algebroids is then easily established.  This furnishes one
answer to the question: When is an arbitrary Lie algebroid a so-called
{\df action} (or {\df transformation}) Lie algebroid?

The remainder of the paper focuses on examples demonstrating the
versatility of the model; in each example the corresponding
implications of Theorem A are identified. This neatly unifies several
results for the first time, although the results themselves are mostly
known.  More substantial applications are pursued in \cite{Blaom_A}.

In Sect.\,\ref{examples} we discuss absolute parallelisms and the
question: When is a differential manifold a Lie group?  More
generally, we turn to the question: When are the leaves of a foliation
the orbits of some Lie group action?  After describing the simplest
scenario explicitly, we sketch how the question might be answered more
generally. Known conditions for the maximal local homogeneity of a
Riemannian manifold are recovered and then generalized in Theorem B
(Sect.\,\ref{section6}) to those infinitesimal $G$-structures
supporting a Cartan connection of `reductive' type.  In
\ref{symplectic} we discuss the existence of {\em finite}-dimensional
symmetries in the Poisson category, again restricting to the simplest
scenario. This gives an answer (admittedly simple-minded) to the
question: When is a Poisson manifold the dual of a Lie algebra?

Sect.\,\ref{section4} reviews the theory of connections from the Lie
algebroid viewpoint and in Sect.\,\ref{section5} we recall how the
first jet bundle associated with a Lie algebroid is another Lie
algebroid; we use this to explain the real meaning of compatibility
mentioned above.  Sect.\,\ref{section5} concludes with a brief
discussion of invariant differential operators associated with a
Cartan algebroid.  This appears to be closely related to the so-called
`tractor' or `local twistor' calculus used in conformal geometry
\cite{Cap_Gover_02}.

In Sect.\,\ref{section7} we explain how Cartan algebroids generalize
Cartan geometries.  We formulate Theorem C,
characterizing those Cartan geometries whose corresponding Cartan
algebroids are locally symmetric, a notion weaker than local flatness
of the Cartan geometry itself (indicating local coincidence with the
prescribed model).

Choosing a Cartan connection on a {\em transitive} Lie algebroid
$\mathfrak g$ amounts to choosing a certain `representation' of
$\mathfrak g$ on itself; see \ref{tr}. In \ref{red} we reduce the
existence of such a representation to the existence of a representation
of $\mathfrak g$ on $TM$. (For infinitesimal $G$-structures, for
example, such a representation is god-given.)  Existence in general is
not addressed.  



The present work concerns {\em infinitesimal} symmetries exclusively.
Lie group symmetries or Lie pseudogroups of symmetries are not discussed.
Efficient tools for globalizing our results would be the groupoid
versions of Lie's Fundamental Theorems \cite{Crainic_Fernandes_03} and
the integrability results of Dazord \cite{Dazord_97}.

Kirill Mackenzie has pointed out that the Kumpera-Spencer theory of
Lie equations \cite{Kumpera_Spencer_72} is a model of geometric
structures based on deformations. A related theory is that of
Griffiths \cite{Griffiths_64}. The relationship between these theories
and the present framework has not been explored.

\subsection*{Background} 
Cartan geometries are not so well-known and their history is somewhat
murky.  Fortunately, a lot may be learned about them from Sharpe's
beautiful book \cite{Sharpe_97}, one inspiration for the present work.
Other recent discussions of Cartan geometries are
\cite{Alekseevsky_Michor_95} and \cite{Slovak_97}. 

Lie algebroids and Lie groupoids are presently the subject of
considerable attention. One highlight is the recent generalization of
Lie's Third Fundamental Theorem, due to Crainic and Fernandes
\cite{Crainic_Fernandes_03}. Our main sources of information have been
\cite{CannasdaSilva_Weinstein_99} and \cite{Mackenzie_03}.
\section{Deforming symmetries}\lab{section2}
\subsection{Infinitesimal symmetries}\lab{infin}
We will be deforming infinitesimal symmetries.  For us these are
infinitesimal actions on a smooth connected manifold $M$ by a Lie
algebra ${\mathfrak g}_0$, i.e., a Lie algebra homomorphism
\begin{equation}
        {\mathfrak g}_0\rightarrow\Gamma(TM)\mathlab{eq1},
\end{equation}
into the space of smooth vector fields on $M$. Only {\em smooth}
actions will be of interest, i.e., those for which the corresponding
action map ${\mathfrak g}_0\times M\xrightarrow{\#}TM$ is smooth.

View ${\mathfrak g}_0\times M$ as a vector bundle over $M$ and use the
same symbol, $\#$, to annotate the associated map of sections,
\begin{equation}
\Gamma({\mathfrak g}_0\times M)\xrightarrow{\#}\Gamma(TM).\mathlab{eqt}  
\end{equation}
This map is an extension of \eqref{eq1} when we regard ${\mathfrak
  g}_0$ as the subspace of $\Gamma({\mathfrak g}_0\times M)$
consisting of constant sections.

A key observation is that the Lie bracket on ${\mathfrak g}_0$ extends
in a natural way to a Lie bracket on $\Gamma ({\mathfrak g}_0\times
M)$ in a way making the extension \eqref{eqt} of \eqref{eq1} a Lie
algebra homomorphism as well.

To obtain the new bracket, one first extends the bracket
$[\,\cdot\,,\,\cdot\,]_{{\mathfrak g}_0}$ on ${\mathfrak g}_0$ in a
trivial way: viewing sections of ${\mathfrak g}_0\times M$ as
${\mathfrak g}_0$-valued functions on $M$, define
$\tau(X,Y)(m):=[X(m),Y(m)]_{{\mathfrak g}_0}$ for sections $X$ and $Y$ of
${\mathfrak g}_0\times M$.  If $\nabla$ denotes the canonical flat
affine connection on ${\mathfrak g}_0\times M$, then the sought after
bracket is defined by
\begin{equation}
        [X,Y]_{{\mathfrak g}_0\times M}:=\nabla_{\#X}Y-\nabla_{\#Y}X 
         + \tau(X,Y),\qquad X, Y\in \Gamma ({\mathfrak g}_0\times M).
         \mathlab{ghj}
\end{equation}
Notice that $\tau$ is just the unique extension of
$[\,\cdot\,,\,\cdot\,]_{{\mathfrak g}_0}$ that is linear with respect
to smooth functions $f$ on $M$. By contrast,
$[\,\cdot\,,\,\cdot\,]_{{\mathfrak g}_0\times M}$ satisfies a Leibnitz
property mimicking that of the Jacobi-Lie bracket on $\Gamma(TM)$:
\begin{equation*}
        [X,fY]_{{\mathfrak g}_0\times M}=f[X,Y]_{{\mathfrak g}_0\times 
        M}+df(\#X)Y.
\end{equation*}
\subsection{Lie algebroids}
Lie algebroids generalize the infinitesimal symmetries described
above.  By definition, a Lie algebroid is any vector bundle
${\mathfrak g}$ over $ M$ (generalizing $\mathfrak g_0\times M$
above), together with a vector bundle morphism $\#\colon{\mathfrak
  g}\rightarrow TM$, called the {\df anchor}, and an ${\mathbb
  R}$-linear Lie bracket $[\,\cdot\,,\,\cdot\,]_{\mathfrak g}$ on
$\Gamma ({\mathfrak g})$ making $\#\colon \Gamma ({\mathfrak
  g})\rightarrow \Gamma (TM)$ into a Lie algebra homomorphism.
Additionally, the bracket should be Leibnitz in the sense that
\begin{equation*}
                [X,fY]_{\mathfrak g}=f[X,Y]_{\mathfrak g}+df(\#X)Y.
\end{equation*}
A {\df morphism} $\phi\colon\mathfrak g\rightarrow\mathfrak h$ of Lie
algebroids is a morphism of the underlying vector bundles (covering the
identity) whose lift $\phi\colon\Gamma(\mathfrak
g)\rightarrow\Gamma(\mathfrak h)$ to sections is a homomorphism of Lie
algebras. Additionally, one requires $\#\circ\phi=\#$.

The bundle ${\mathfrak g}_0\times M$ described in \ref{infin} is
called an {\df action Lie algebroid}. In analogy with this case, an
arbitrary Lie algebroid ${\mathfrak g}$ is {\df transitive} if its
anchor $\#\colon\mathfrak g\rightarrow TM$ is surjective; if $\mathfrak
g$ is intransitive, the (possibly singular) distribution
$\#({\mathfrak g})\subset TM$ is nevertheless integrable, giving rise
to a foliation by the {\df orbits} of the algebroid. A Lie algebroid
$\mathfrak g$ is {\df regular} if the orbits have constant dimension,
or equivalently, if the anchor has a subbundle of $\mathfrak g$
as kernel.

The tangent bundle $TM$ is itself a Lie algebroid over $M$, with the
identity on $TM$ as anchor; every Lie {\em algebra} is a Lie algebroid
over a single point. A Lie algebroid with trivial anchor is a bundle
of Lie algebras. In particular, this applies to the kernel of the
anchor of a regular Lie algebroid. See Sect.\,\ref{examples} and work
cited above for further examples.

\subsection{Cartan algebroids}\lab{geom}
The connection $\nabla$ appearing in \ref{infin} is an instance of the
following general notion, central to all that follows:
\begin{definition}
  A {\df Cartan connection} on a Lie algebroid ${\mathfrak g}$ over
  $M$, with anchor $\#$, will be any affine connection
  $\nabla$ on $\mathfrak g$ that is compatible with its bracket
  $[\,\cdot\,,\,\cdot\,]_{{\mathfrak g}}$ in the following sense:
  \begin{equation}
    \nabla_V[X,Y]_{\mathfrak g}=[\nabla_VX,Y]_{\mathfrak g}
                              +[X,\nabla_VY]_{\mathfrak g}
                              +\nabla_{\bar\nabla_YV}X
                              -\nabla_{\bar\nabla_XV}Y,\mathlab{compat}
  \end{equation}
for all $X,Y\in\Gamma({\mathfrak g})$ and $V\in\Gamma(TM)$. Here
\begin{equation*}
  \bar\nabla_XV:=\#\nabla_VX+[\#X,V]_{TM},\qquad X\in\Gamma(\mathfrak g),
                                                 V\in\Gamma(TM),
\end{equation*}
where $[\,\cdot\,,\,\cdot\,]_{TM}$ denotes the Jacobi-Lie bracket on
vector fields on $M$. A {\df Cartan algebroid} is a Lie algebroid
equipped with a Cartan connection. A {\df morphism} of Cartan
algebroids is simply a connection-preserving
morphism of the underlying Lie algebroids.
\end{definition}

Regarding \eqref{compat}: An affine connection on $\mathfrak g$
amounts to a splitting $J^1\mathfrak g\xleftarrow{s}\mathfrak g$ of
a canonical exact sequence 
$    0\rightarrow T^*M\otimes\mathfrak g\monomorphism J^1\mathfrak g
    \rightarrow\mathfrak g\rightarrow 0.$ 
Here $J^1\mathfrak g$ denotes the first jet bundle of $\mathfrak g$.
As we elaborate in Sect.\,\ref{section5}, each arrow is a morphism of
Lie algebroids and condition \eqref{compat} holds precisely when
$J^1\mathfrak g\xleftarrow{s}\mathfrak g$ is a Lie algebroid morphism
also.
  
The operator $\bar\nabla$ is an example of {\df $\mathfrak
  g$-connection} on $TM$; see \ref{gConn}.
The $\mathfrak g$-connection is {\df flat} if
\begin{equation*}
      (\bar\nabla_X\bar\nabla_Y-\bar\nabla_Y\bar\nabla_X
      -\bar\nabla_{[X,Y]_{\mathfrak g}})V=0,
\end{equation*}
for all $X,Y\in\Gamma(\mathfrak g)$ and $V\in\Gamma(TM)$. Flat
$\mathfrak g$-connections are called {\df $\mathfrak
  g$-representations} because they generalize the usual
representations of a Lie algebra; instead of acting on a vector space,
Lie algebroids act on vector {\em bundles}, in this case $TM$. 

If $\nabla$ is to be a Cartan connection then it is necessary that
$\bar\nabla$ be a $\mathfrak g$-representation. One sees this by
applying the anchor $\#$ to both sides of \eqref{compat}.  If
$\#\colon\mathfrak g\rightarrow TM$ is injective, this condition is
sufficient.  (For examples, see \ref{abs} and \ref{foliations}.) Lie
algebroid representations are reviewed in Sect.\,\ref{section4}.

From a Cartan connection $\nabla$ one also obtains a representation 
of $\mathfrak g$ on itself, with respect to which the anchor $\#$ is 
equivariant.  See \ref{pat}.  This is significant because, unlike Lie 
algebras, a Lie algebroid does not generally admit an adjoint 
representation in the naive sense (but do see \ref{adjoint}).

\subsection{The symmetric part of a Cartan algebroid}\lab{substructure}
A fundamental observation is that every Cartan algebroid $(\mathfrak
g,\nabla)$ has a canonical subalgebroid isomorphic to an action Lie
algebroid. Indeed, let ${\mathfrak g}_0\subset\Gamma(\mathfrak g)$ be
the subspace of $\nabla$-parallel sections, which is
finite-dimensional. Then bracket compatibility \eqrefs{geom}{compat}
guarantees that ${\mathfrak g}_0\subset\Gamma(\mathfrak g)$ is a Lie
subalgebra, and we obtain an action of $\mathfrak g_0$ on $M$ given by
\begin{gather*}
  \mathfrak g_0\times M\rightarrow TM\\
  (X,m)\mapsto \# X(m).  
\end{gather*}
Equipping $\mathfrak g_0\times M$ with the canonical flat affine connection,
we obtain a morphism of Cartan algebroids,
\begin{gather}
  \mathfrak g_0\times M\rightarrow\mathfrak g\label{kernel}\\
  (X,m)\mapsto X(m).\notag
\end{gather}
This morphism is injective because $\nabla$-parallel sections
vanishing at a point vanish everywhere. (We are assuming $M$ is
connected.) We call the image of the monomorphism \eqref{kernel} the
{\df symmetric part} of $(\mathfrak g, \nabla)$.

\subsection{Curvature as the local obstruction to symmetry}\lab{dissymmetry}%
A Cartan algebroid $(\mathfrak g, \nabla)$ on $M$ is {\df
  symmetric} if it is isomorphic to an action algebroid ${\mathfrak
  g}_0\times M$, equipped with its canonical flat connection --- or
equivalently, if it coincides with its symmetric part. We call
$(\mathfrak g,\nabla)$ {\df locally symmetric} if every point of $M$
has an open neighborhood $U$ on which the restriction $(\mathfrak g,
\nabla)|_U$ is symmetric. 

Cartan algebroids are indeed symmetries deformed by curvature:
\begin{theoremA}\lab{thmA}
  A Cartan algebroid $(\mathfrak g, \nabla)$ on $M$ is locally
  symmetric if and only if $\nabla$ is flat. When $M$ is simply
  connected, local symmetry already implies symmetry.
\end{theoremA}
\begin{proof}
  The necessity of flatness is immediate.  To finish the proof it 
  suffices to show that \eqref{kernel} is an isomorphism whenever 
  $\nabla$ is flat and $M$ is simply connected.  Indeed, in that case 
  $\nabla$ determines a trivialization of the bundle $\mathfrak g$ in 
  which constant sections correspond to the $\nabla$-parallel sections 
  of $\mathfrak g$ --- that is, to elements of $\mathfrak g_0$.  (This 
  classical result follows, for example, from the groupoid version of 
  Lie's Second Fundamental Theorem; see Remark \eqrefs{TM}{TM2}.)  In 
  particular, $\mathfrak g_0\times M$ and $\mathfrak g$ will have the 
  same rank, implying the monomorphism \eqref{kernel} is an 
  isomorphism.
\end{proof}

\section{Examples}\lab{examples}
We now describe simple examples of Cartan algebroids and describe
simple consequences of Theorem A, stated as corollaries.  More general
classes of examples are presented in Sections \ref{section6} and
\ref{section7}.

\subsection{Absolute parallelisms}\lab{abs}
The simplest application of Theorem A is to the question: When is a
manifold a Lie group?

Affine connections on a tangent bundle $TM$ occur in pairs.  If
$\nabla$ is one such connection, its {\df dual} $\nabla^*$ is defined
by
\begin{equation*}
        \nabla_X^*Y:=\nabla_YX+[X,Y]_{TM}.
\end{equation*}
We have $\nabla^{**}=\nabla$ and call $ (\nabla,\nabla^*)$ a {\df dual
  pair}.  When $M$ is a Lie group, the flat connections associated
with the canonical left and right trivializations of the tangent
bundle $TM$ are a dual pair in this sense. Conversely, any dual pair
of simultaneously flat connections determines a local Lie group
structure, as we now explain. 

A Cartan connection on $TM$ is simply an affine connection whose dual
is flat. This follows from comments made in \ref{geom} (where
$\nabla^*=\bar\nabla$).  So specifying a Cartan connection is
equivalent to specifying a flat affine connection on $TM$.  These
latter connections are the infinitesimal analogues of absolute
parallelisms (trivializations of the tangent bundle), the most basic
of all geometric structures.

Specifying an absolute parallelism on $M$ is equivalent to specifying
a one-form $\omega$ on $M$, taking values in some vector space $ V$,
such that $\omega\colon T_mM\rightarrow V$ is an isomorphism at all
$m\in M$. 
\begin{corollary}
  Let $(\nabla,\nabla^*)$ be a dual pair of connections on the tangent
  bundle of a simply connected manifold $M$ and suppose $\nabla$ and
  $\nabla^*$ are both flat. Then there exists a Lie algebra
  ${\mathfrak g}_0$ and an absolute parallelism
  $\omega\in\Omega^1(M,{\mathfrak g}_0)$ such that
        \begin{gather*}
                 d\omega+\frac{1}{2}[\omega,\omega]_{{\mathfrak g}_0}=0,\\
                 \omega(\nabla_XY)=(\mathcal L_X\omega)(Y)\quad
              \text{and}\quad\omega(\nabla^*_XY)=\mathcal 
              L_X(\omega(Y)).
        \end{gather*}
  Here $\mathcal L$ denotes Lie derivative.
\end{corollary}
\noindent If in addition the `Mauer-Cartan' form $\omega$ is {\em complete},
then $M$ is diffeomorphic to the simply connected Lie group $ G_0$
integrating ${\mathfrak g}_0$; see \cite[Theorem {\bfseries
  3}.8.7]{Sharpe_97} or \cite{Gardner_89}.
\begin{remark}
  If $\nabla^*$ is flat then $\nabla$ is also flat if and only if the
  torsion of $\nabla^*$ is $\nabla^*$-parallel (Proposition
  \eqrefs{dual}{scorch}).
\end{remark}

The corollary is established by taking $\mathfrak g=TM$ in Theorem A,
which then delivers an isomorphism $TM\rightarrow\mathfrak g_0\times
M$.  One takes $\omega$ to be the composite
$TM\rightarrow\mathfrak g_0\times M\rightarrow\mathfrak g_0$.

\subsection{Foliations}\lab{foliations}
It is natural to ask when the leaves of an arbitrary foliation
$\mathcal F$ are the orbits of some action by a Lie group, or at least
a Lie algebra. We now mention the most elementary application of
Theorem A to this question.

Let $D\subset TM$ be a subbundle that is integrable as a distribution.
Then $D\subset TM$ is a Lie subalgebroid. An affine connection
$\nabla$ on $D$ is a Cartan connection precisely when the
corresponding $D$-connection $\bar\nabla$, defined by
\begin{equation*}
  \bar\nabla_XV:=\nabla_VX+[X,V]_{TM},\qquad X\in\Gamma(D),V\in\Gamma(TM),
\end{equation*}
is flat. Call a Lie algebra action {\em free} if all orbits have the
dimension of the Lie algebra. Theorem A implies:
\begin{corollary}
  Let $\mathcal F$ be a regular foliation on a simply connected
  manifold $M$, and let $D\subset TM$ denote its tangent distribution.
  Then the orbits of $\mathcal F$ are the orbits of a free Lie algebra
  action if and only if $D$ admits a flat affine connection $\nabla$
  whose corresponding $D$-connection $\bar\nabla$ is also flat.
\end{corollary}

Taking $ D=TM$, we recover the preceding corollary, since the
existence of a Mauer-Cartan form is nothing more than the existence of
a free and transitive action by a Lie algebra.  

More generally, one may look for Lie algebra actions of `higher order'
and consider smooth foliations ${\mathcal F}$ with singularities.
Call an action by a Lie algebra ${\mathfrak g}_0$ {\df $k$th-order
  faithful} if the natural morphism ${\mathfrak g}_0\times
M\rightarrow J^k(TM)$ into the $k$th-order jet bundle of $TM$ is
injective, i.e., realizes $\mathfrak g_0\times M$ as a subalgebroid of
$J^k(TM)$.  If $ D\subset TM$ is the distribution tangent to a
possibly singular foliation ${\mathcal F}$, then one searches for
subalgebroids ${\mathfrak g}\subset J^k (TM)$ with $\#({\mathfrak
  g})=D$, and establishes `$k$th-order symmetry' of $\mathcal F$ by
finding a flat Cartan connection on ${\mathfrak g}$.  In the case
$D=TM$, obstructions in this search represent obstructions to the
realization of $M$ as a homogenous space. 

\subsection{Riemannian structures}\lab{Riemannian}
Let $\sigma$ be a Riemannian metric on $M$. Then an
infinitesimal isometry of $(M, \sigma)$, or {\df Killing field}, is
a vector field $V$ on $M$ such that $\mathcal L_V \sigma=0$, where $\mathcal
L$ is Lie derivative. That is,
\begin{equation*}
  (\mathcal L_V \sigma)(U,W):= 
  \mathcal L_V(\sigma(U,W))-\sigma(\mathcal L_VU,W)+ \sigma(U,\mathcal L_VW).
\end{equation*}
Generally speaking, no Killing fields exist, even locally.  However,
for each fixed $m\in M$, there are many vector fields $V$ such that
$\mathcal L_V \sigma$ vanishes {\em at the point $m$.} The one-jets
$[V]_m^1$ at $m$ of all such $V$ are elements of a fundamental Lie
algebroid $\mathfrak g$ associated with $ \sigma$. In detail, if
$J^1(TM)$ denotes the first jet bundle of $TM$ (a Lie algebroid, see
\ref{structure}), then we define $\mathfrak g$ to be the kernel of
\begin{gather*}
  J^1(TM)\rightarrow\symmetric^2(TM)\\
  [V]_m^1\mapsto (\mathcal L_V \sigma)(m).
\end{gather*}
This is a transitive Lie subalgebroid of $J^1(TM)$. The anchor $\#$ of
$\mathfrak g$ is the restriction of the natural projection
$J^1(TM)\rightarrow TM$.  The kernel $\mathfrak h\subset T^*M\otimes
TM\subset J^1(TM)$ of this anchor is the Lie algebra bundle of all
$\sigma$-skew symmetric tangent space endomorphisms.

Because $\mathfrak g\subset J^1(TM)$ is transitive, there exists a
(non-unique) splitting $\mathfrak g\xleftarrow{t}TM$ of the exact
sequence $0\rightarrow\mathfrak h \monomorphism \mathfrak
g\stackrel{\#}{\rightarrow}TM\rightarrow 0$ such that $\mathfrak
g=\mathfrak h+t(TM)$. For every such $t$, the composite
$J^1(TM)\hookleftarrow\mathfrak g\xleftarrow{t}TM$ is a splitting of
the exact sequence $0\rightarrow T^*M\otimes TM\monomorphism
J^1(TM)\rightarrow TM\rightarrow 0$, and so determines an affine
connection $\levi$ on $TM$. This is precisely the Levi-Cevita
connection corresponding to $\sigma$ when $t$ is the unique splitting
making $\levi$ torsion free.

As we explain in \ref{adjoint}, there is a natural representation of
$J^1(TM)$ on $TM$, denoted $\adjoint$, and defined by $\adjoint_{J^1
  V}W=[V,W]_{TM}$, where $ J^1 V$ denotes the first prolongation of $
V$.  This representation restricts to a representation of $\mathfrak
g$ on $TM$.  As we detail in \ref{red}, once a representation of
$\mathfrak g$ on $TM$ is prescribed, a splitting $t$ as above 
determines an associated `reductive' Cartan connection
$\nabla$ on $\mathfrak g$. Choosing the splitting $t$ corresponding to
the Levi-Cevita connection, we obtain a connection $\nabla$ with the
following property:
\begin{proposition}
  A vector field $V$ is a Killing field for $\sigma$ if and only if
  $V=\#X$ for some $\nabla$-parallel section $X$ of $\mathfrak g$
  (equivalently, if $J^1V$ is a section of $\mathfrak g$).
\end{proposition}
\noindent%
Sufficiency of the stated condition is not difficult to see. Its
necessity may be established using Cartan's method of equivalence. See
\cite{Blaom_A}, which contains an alternative construction of $\nabla$
(there denoted $\nabla^{(1)}$).

The proposition establishes a correspondence between the
Killing fields of $\sigma$ and the $\nabla$-parallel sections of
$\mathfrak g$. It implies that $(M,\sigma)$ is locally maximally
homogeneous, in the sense of locally possessing a Lie algebra of
Killing fields of maximal dimension, if and only if $(\mathfrak
g,\nabla)$ is locally symmetric, in the sense of \ref{dissymmetry}.
\begin{corollary}
  The above Cartan algebroid $(\mathfrak g,\nabla)$ on a Riemannian
  manifold $(M,\sigma)$ is locally maximally homogeneous, if and
  only if the curvature of the Levi-Cevita connection is
  simultaneously $\mathfrak h$-invariant, and $\levi$-parallel.
\end{corollary}
\noindent%
As is
well-known, the stated condition is equivalent to pure, constant scalar
curvature, or equivalently, to constant sectional curvature.
This corollary of Theorem A is a special case of Theorem B,
p.\,\pageref{here}.  

\subsection{Symplectic and Poisson structures}\lab{symplectic}
Let $\omega$ be a symplectic structure on $M$ and let $\#\colon
T^*M\rightarrow TM$ denote the inverse of
$v\mapsto\omega(v,\,\cdot\,)$.  Since $\#$ is an isomorphism, there is
a unique bracket on $\Gamma(T^*M)$ making $T^*M$ into a Lie algebroid
with anchor $\#$. This bracket is given by
\begin{equation}
  [\alpha,\beta]_{T^*M}=\mathcal L_{\#\alpha}\beta-\mathcal L_{\#\beta}\alpha 
                       + d(\Pi(\alpha,\beta)),\qquad\alpha,\beta\in\Gamma(T^*M),
  \mathlab{one}
\end{equation}
where $\mathcal L$ denotes Lie derivative and $\Pi$ is the Poisson
tensor. This tensor is defined by
$\Pi(\alpha,\beta):=\omega(\#\alpha,\#\beta)$ and so satisfies
\begin{equation}
  \langle\alpha,\#\beta\rangle=\Pi(\alpha,\beta)
  \qquad\alpha,\beta\in\Gamma(T^*M).\mathlab{two}
\end{equation}
As is well-known, \eqref{one} defines a Lie algebroid structure on
$T^*M$ for {\em any} Poisson manifold $(M,\Pi)$, with anchor $\#$ defined by
\eqref{two}.

An  infinitesimal isometry of a Poisson manifold $(M,\Pi)$ is a
vector field $V$ on $M$ such that $\mathcal L_V\Pi=0$. Like
foliations, Poisson manifolds have an abundance of infinitesimal
isometries. In particular, every {\em closed} one-form $\alpha$ on $M$
determines an infinitesimal symmetry $\#\alpha$ tangent to the
symplectic leaves known as a {\df local Hamiltonian vector
  field}, or a {\df Hamiltonian vector field} if $\alpha$ is exact.

The simplest instance of finite-dimensional symmetry occurs when
$T^*M$ is itself isomorphic to an action Lie algebroid. This occurs,
for instance, when $M$ is the dual of a Lie algebra, equipped with its
Lie-Poisson structure \cite[\S 10.1]{Marsden_Ratiu_94}.  According to
the corollary below, this is almost the only such case.

Let $\nabla$ be an affine connection on the Lie algebroid $T^*M$
associated with an arbitrary Poisson manifold $(M,\Pi)$. Assume that
the space $\mathfrak g_0$ of $\nabla$-parallel sections consists of
closed one-forms. It suffices to suppose $\nabla$ is {\df torsion
  free}, in the sense that the corresponding connection on $TM$ is
torsion free; in the flat case this sufficient condition becomes
necessary as well.  Equivalently, 
  \begin{equation}
    d\alpha(V,W)=\langle\nabla_V\alpha,W\rangle
    -\langle\nabla_W\alpha,V\rangle,\mathlab{p1}    
  \end{equation}
for all one-forms $\alpha$ and vector fields $V,W$ on $M$.
  
Assuming $\nabla$ is torsion free, we can readily use \eqref{p1} to
rewrite \eqref{one} as
\begin{equation}
  [\alpha,\beta]_{T^*M}=\nabla_{\#\alpha}\beta-\nabla_{\#\beta}\alpha
  -\nabla\Pi\,(\alpha,\beta),\mathlab{p2}
\end{equation}
where $\nabla\Pi$ is
the section of $\alternating^2(T^*M,T^*M)$ defined by
\begin{equation*}
  \nabla\Pi\,(\alpha,\beta)V=(\nabla_V\Pi)(\alpha,\beta),
  \qquad\alpha,\beta\in\Gamma(T^*M),\,V\in\Gamma(TM).
\end{equation*}
Using \eqref{p2}, it is not too difficult to characterize those
torsion free connections satisfying \eqrefs{geom}{compat}:
\begin{lemma}
  A torsion free connection $\nabla$ on $T^*M$ is a Cartan
  connection if and only if
  \begin{equation}
    \curvature\nabla\,(V,\#\alpha)\beta-\curvature\nabla(V,\#\beta)\alpha
    -(\nabla_V\nabla\Pi)(\alpha,\beta)=0,\mathlab{sx}
  \end{equation}
  for all one-forms $\alpha,\beta$ and vector fields $V$ on $M$.
\end{lemma}
\noindent%

We may now prove the following corollary of Theorem A, answering the
question: When is a Poisson manifold the dual of a Lie algebra?  Note
that hypotheses \eqref{qq1} and \eqref{qq2} of the corollary are
necessary in this case.
\begin{corollary}
  Let $(M,\Pi)$ be a simply connected Poisson manifold. Assume:
  \begin{conditions}
  \item\lab{qq1} The tangent bundle $TM$ admits a flat, torsion free,
    connection $\nabla$ such that $\nabla\Pi$ is
    $\nabla$-parallel.
  \end{conditions}
  Then the space $\mathfrak g_0$ of $\nabla$-parallel one-forms on $M$
  is a finite-dimensional subalgebra of $\Gamma(T^*M)$, acting
  leaf-transitively on $M$ by Hamiltonian vector fields.  Moreover
  assume:
  \begin{conditions}
  \item\lab{qq2} Some corresponding momentum map $\mathbf J\colon
    M\rightarrow \mathfrak g_0^*$ is proper and infinitesimally
    equivariant $($see, e.g., \cite{Marsden_Ratiu_94}$)$.
  \end{conditions}
  Then $\mathbf J$ maps $M$ isomorphically onto an open subset of
  $\mathfrak g_0^*$, equipped with its Lie-Poisson structure.
\end{corollary}
%
\begin{proof}
  Assuming \eqref{qq1} holds, denote the corresponding affine
  connection on $ T^*M$ by $\nabla$ also. By the lemma, $\nabla$ is a
  Cartan connection on $ T^*M$. That ${\mathfrak g}_0\subset \Gamma
  (T^*M)$ is a subalgebra acting on $ M$ follows from observations
  made in \ref{substructure}. Since $\nabla$ is torsion free and $M$
  is simply connected, elements of ${\mathfrak g}_0$ are exact
  one-forms, so that ${\mathfrak g}_0$ acts via Hamiltonian vector
  fields. The action consequently admits a momentum map.
  
  Since $\nabla$ is flat, the canonical inclusion ${\mathfrak
    g}_0\times M\monomorphism T^*M$ is an isomorphism (see the proof
  of Theorem A), implying that ${\mathfrak g}_0$ acts
  leaf-transitively and that any momentum map ${\mathbf J}\colon
  M\rightarrow \mathfrak g_0^*$ is a local diffeomorphism.  Properness
  of ${\mathbf J}$ and simple connectivity of $M$ then imply ${\mathbf
    J}$ is a diffeomorphism onto its image.  That ${\mathbf J}$ is
  additionally Poisson follows from its infinitesimal equivariance
  \cite[\S 12.4]{Marsden_Ratiu_94}.
\end{proof}
\section{Connections and Lie algebroid representations}\lab{section4}
The next two sections expand on remarks made in \ref{geom} and prepare
us for the general classes of examples to be described in Sections
\ref{section6} and \ref{section7}. The present section reviews basic
connection theory from the Lie algebroid point of view.  

Lie algebroid formalism suggests a generalization of affine
connections called {\df $\mathfrak g$-connections}, studied
systematically in \cite{Fernandes_02} (where they are the {\df
  $A$-derivatives} associated with {\df $A$-connections}). Just as Lie
algebras act on vector spaces, so Lie {\em algebroids} act on vector
{\em bundles}. A $\mathfrak g$-connection is just a Lie algebroid
representation `deformed by curvature' as we presently explain.

What we will describe as the `general linear Lie algebroid $\gl(E)$'
of a vector bundle $E$ has been previously described in
\cite{Higgins_Mackenzie_90}, where it is denoted
$\operatorname{DO}(E)$.

\subsection{Curvature}\lab{curv}
All the `curvatures' in this paper are instances of a single Lie-algebraic
notion. Recall that a vector bundle morphism $\phi\colon\mathfrak
g\rightarrow\mathfrak h$ between Lie algebroids is a Lie algebroid
morphism if:
\begin{conditions}
  \item\lab{redone} $\#\circ\phi=\#$~,\enspace and
  \item\lab{redtwo} $\phi[X,Y]_{\mathfrak g}=[\phi X,\phi Y]_{\mathfrak h}$ for
    all $X,Y\in\Gamma(\mathfrak g)$.
\end{conditions}
When only \eqref{redone} holds, we define the {\df curvature} of $\phi$ by
\begin{equation*}
  \curvature\phi\,(X,Y):=[\phi X,\phi Y]_{\mathfrak h}-\phi[X,Y]_{\mathfrak g}
  \qquad X,Y\in\Gamma(\mathfrak g).
\end{equation*}
Thus $\phi$ is a Lie algebroid morphism if and only if
$\curvature\phi=0$. Notice that $\curvature\phi\,(X,Y)$ always lies in
the kernel of $\#\colon\mathfrak h\rightarrow TM$.

\subsection{The general linear Lie algebroid of a vector 
bundle}\lab{glE}%
A representation of a Lie {\em algebra} ${\mathfrak g}$ is a vector 
space $E$, together with a Lie algebra homomorphism
\begin{equation*}
        {\mathfrak g}\rightarrow\gl (E).
\end{equation*}
Here $\gl( E):=\homomorphism (E, E)$ is the Lie algebra of the general
linear group $\operatorname{GL}(E)=\automorphism (E)$.  If instead $
E$ is a vector {\em bundle} over $M$, then we use
$\operatorname{GL}(E)$ to denote the {\df frame groupoid} of $ E$,
consisting of all isomorphisms $ E_{m_1}\rightarrow E_{m_2}$ between
fibers of $ E$ over possibly distinct points $ m_1, m_2\in M$. We
define $\automorphism (E)\subset\operatorname{GL}(E)$ to be the
subgroupoid consisting of isomorphisms $ E_{m_1}\rightarrow E_{m_2}$
with $ m_1=m_2$ (the isotropy subgroupoid).  This subgroupoid is a
bundle of Lie groups.

We now describe, in purely algebraic terms, a concrete model for the
Lie algebroid of $\operatorname{GL}(E)$, a model we will denote by
$\gl(E)$. Sections of $\gl(E)$ will be characterized as differential
operators.  The Lie algebroid bracket of two sections will simply be
the commutator of the corresponding operators.

Let $ J^1 E$ denote the first jet bundle of $ E$ and let $ J^1\colon
\Gamma (E)\rightarrow \Gamma (J^1 E)$ denote prolongation:
$J^1(\sigma)(m):=[\sigma]_m^1$. Here $[\sigma]_m^1$ denotes the
one-jet of $\sigma$ at $m$.  We have an exact sequence
\begin{equation}
   0\rightarrow T^*M\otimes E\monomorphism 
   J^1E\rightarrow E\rightarrow 0,\mathlab{crazy}
\end{equation}
where the inclusion on the left is the one sending a section
$df\otimes \sigma$ of $ T^*M\otimes E$ to the section $ J^1( f\sigma)-f J^1 \sigma$
of $J^1 E$. Applying $\homomorphism (\,\cdot\,, E)$ to the sequence,
and identifying $\homomorphism (T^*M\otimes E, E)$ with
$TM\otimes\homomorphism (E, E)$, we obtain a second exact sequence
\begin{equation*}
   0\rightarrow\homomorphism (E, E)\monomorphism\homomorphism (J^1 E,E)
   \xrightarrow{*}TM\otimes\homomorphism(E, E)\rightarrow 0.
\end{equation*}
\begin{definition}
  The {\df general linear Lie algebroid} $\gl(E)$ of a vector bundle
  $E$ over $M$ is the preimage, under the surjective arrow $*$, of
  the subbundle $TM\otimes\{\identity_E\}\subset
  TM\otimes\homomorphism(E,E)$.
\end{definition}
The general linear Lie algebroid $\gl(E)$ of $E$ is indeed a Lie algebroid, 
as we now explain.  Firstly, its anchor is the map $\#$ completing the 
following commutative diagram:
\begin{equation*}
        \begin{CD}
                \gl(E)@>\#>> TM\\
                @V\text{inclusion}VV @VVv\mapsto v\otimes\identity_EV\\
                \homomorphism (J^1 E,E)@>*>>TM\otimes\homomorphism(E,E)
        \end{CD}
\end{equation*}
The bundle of Lie algebras $\homomorphism (E, E)\subset\homomorphism
(J^1 E, E)$ (the Lie algebroid of $\automorphism(E)$) is contained in
$\gl(E)$.  Moreover $\homomorphism (E, E)$ is precisely the kernel of
the anchor $\#$ above, i.e., we have an exact sequence
\begin{equation*}
        0\rightarrow\homomorphism (E, E)\monomorphism\gl 
        (E)\rightarrow{\#}TM\rightarrow 0.
\end{equation*}

For any section $\xi$ of $\homomorphism (J^1 E, E)$, we define a first 
order differential operator, $\operator (\xi)\colon \Gamma 
(E)\rightarrow \Gamma (E)$, by
\begin{equation*}
        \operator (\xi) \sigma:=\xi (J^1 \sigma).
\end{equation*}
Unravelling the definitions above, one has:
\begin{proposition}
  A section $\xi$ of $\homomorphism(J^1E,E)$ is a section of $\gl(E)$
  if and only if there exists a vector field $ V$ on $M$ such that
         \begin{equation*}
                \operator (\xi) (f\sigma)=f\operator (\xi) \sigma+df(V)\sigma
         \end{equation*}
         for all smooth functions $ f$ on $M$, and sections $\sigma$ 
         of $E$.
         In that case, $V=\#\xi$.
\end{proposition}

A Lie algebroid bracket for $\gl(E)$ is now well defined by
\begin{equation*}
        \operator ([\xi,\eta]_{\gl(E)})=\operator(\xi)\circ\operator 
        (\eta)-\operator(\eta)\circ\operator(\xi).
\end{equation*}
By construction, the vector bundle $\homomorphism (E, E)$ is contained 
in $\gl(E)$ as a Lie subalgebroid.  The bracket 
$[\,\cdot\,,\,\cdot\,]_{\homomorphism (E,E)}$ that $\homomorphism 
(E,E)$ inherits from that on $\gl(E)$, agrees with its usual Lie algebra 
bundle bracket:
\begin{equation}
       [\phi_1,\phi_2]_{\homomorphism (E, E)}
       =\phi_1\circ\phi_2-\phi_2\circ\phi_1.
       \mathlab{cand}
\end{equation}

It can be shown that the Lie algebroid $\gl(E)$ is indeed a model for
the (abstract) Lie algebroid of $\operatorname{GL}(E)$.  An element of
the latter may be regarded as some tangent vector ${\dot g}(0)\in
T\operatorname{GL}(E)$, where the path $t\mapsto
g(t)\in\operatorname{GL}(E)$ consists of isomorphisms $ g(t)\colon
E_{m(t)}\rightarrow E_m$ with $ m(0)=m$ (i.e., is an `infinitesimal moving
$E$-frame').  The corresponding element of
$\gl(E)\subset\homomorphism (J^1 E, E)$ sends $[\sigma]_m^1\in J^1 E$ to
$\frac{d}{dt}g(t)\sigma(m(t)){\vert_{t=0}}\in E$.

\subsection{${\mathfrak g}$-representations}\lab{reps}%
A {\df representation} of a Lie algebroid ${\mathfrak g}$ is a vector
bundle $ E$, together with a morphism of Lie algebroids, ${\mathfrak
  g}\rightarrow\gl(E)$.  If ${\mathfrak g}$ and $E$ are vector bundles
over a single point, we obtain the usual notion of a Lie algebra
representation.  For another example, consider a principal bundle
$P\rightarrow M$ with a simply connected structure group $H$, and
corresponding Lie algebroid ${\mathfrak g}:=(TP)/H$, described in
\ref{upi}.  Then the representations of ${\mathfrak g}$ are the
associated vector bundles of $ P$; see, e.g., \cite{Mackenzie_03}.  In
particular, there is in this case a one-to-one correspondence between
(isomorphism classes of) representations of ${\mathfrak g}$ and vector
space representations of $ H$.

\subsection{$\mathfrak g$-connections}\lab{gConn}%
Given a Lie algebroid ${\mathfrak g}$, a {\df ${\mathfrak
    g}$-connection} on a vector bundle $ E$ is a morphism of vector
bundles
\begin{equation*}
        \nabla\colon{\mathfrak g}\rightarrow\gl(E)
\end{equation*}
satisfying $\#\circ\nabla=\#$ that is possibly {\em not} a Lie algebroid
morphism.  Recall that $\gl(E)\subset\homomorphism(J^1E,E)$. For any
vector bundle morphism
\begin{equation*}
  \nabla\colon\mathfrak g\rightarrow\homomorphism(J^1E,E)  
\end{equation*}
one writes $\nabla_X\sigma:=\nabla(X)(J^1\sigma)$
($X\in\Gamma(\mathfrak g)$, $\sigma\in\Gamma(E)$). The following
elementary result characterizes $\mathfrak g$-connections as certain
differential operators, and shows that $\mathfrak g$-representations
are just flat $\mathfrak g$-connections, where flatness takes a
familiar form.
\begin{proposition}\mbox{}
  \begin{conditions}
  \item $\nabla\colon\mathfrak 
    g\rightarrow\homomorphism(J^1E,E)$ is a $\mathfrak g$-connection if 
    and only if
    \begin{equation*}
      \nabla_X(f\sigma)=f\nabla_X\sigma+df(\#X)\sigma
    \end{equation*}
    for all smooth functions $f$ on $M$ and sections
    $X\in\Gamma(\mathfrak g)$, $\sigma\in\Gamma(E)$.
    
  \item In that case $\nabla\colon\mathfrak g\rightarrow\gl(E)$ is 
    a $\mathfrak g$-representation if and only if its curvature,
    in the sense of \ref{curv}, vanishes. This curvature is a section
    of $\alternating^2(\mathfrak g, \homomorphism (E, E))$ and
    satisfies
    \begin{equation*}
        \curvature\nabla\, (X,Y)\sigma
         =(\nabla_X\nabla_Y-\nabla_Y\nabla_X-\nabla_{[X,Y]}) \sigma.
    \end{equation*}
  \end{conditions}
\end{proposition}

\subsection{$TM$-connections}\lab{TM}
A $ TM$-connection is simply an affine connection in the usual sense.
Given a $TM$-connection $\nabla\colon TM\rightarrow\gl( E)$ on $E$,
and a section $\sigma$ of $ E$, we write $\nabladot\sigma$ for the
section of $T^*M\otimes E\cong\homomorphism (TM, E)$ defined by
\begin{equation*}
         (\nabladot\sigma) V:=\nabla_V\sigma=\nabla (V)(J^1\sigma),\qquad
         V\in\Gamma(TM).
\end{equation*}
Then $\nabla$ determines a 
splitting $ J^1E\xleftarrow{s} E$ of the exact sequence 
\eqrefs{glE}{crazy}, given by
\begin{equation}
         s\sigma=J^1\sigma-\nabladot\sigma,
          \qquad\sigma\in\Gamma(E).\mathlab{cost}
\end{equation}
Conversely, \eqref{cost} implicitly defines a $TM$-connection $\nabla$ 
for every such splitting $s$.
\begin{remarks}\mbox{}
        \begin{conditions}
        \item\lab{TM1} The preceding observations are transparent if
          one observes that prolongation $ \Gamma
          (J^1E)\xleftarrow{J^1} \Gamma (E)$ is a god-given splitting
          for the sequence obtained by applying $\Gamma(\,\cdot\,)$ to
          \eqrefs{glE}{crazy}:
         \begin{equation*}
                 0\rightarrow\Gamma( T^*M\otimes E)\monomorphism 
                 \Gamma (J^1 E)\rightarrow
                 \Gamma (E)\rightarrow 0.
         \end{equation*}
                 
       \item\lab{TM2} Suppose $M$ is simply connected.  Then a
         representation of $ TM$ is just a vector bundle $ E$,
         together with a trivialization $ E\cong E_0\times M$.
         Indeed, suppose $\nabla\colon TM\rightarrow\gl(E)$ is such a
         representation (i.e., a flat $ TM$-connection on $ E$).  By
         the groupoid version of Lie's Second Fundamental Theorem
         (see, e.g., \cite{Crainic_Fernandes_03}), the Lie algebroid
         morphism $\nabla$ lifts to a groupoid morphism $\Phi\colon
         M\times M\rightarrow\operatorname{GL}(E)$.  If $ m_0\in M$ is
         fixed and $ E_0$ is the fiber over $ m_0$, the trivialization
         is given by
                \begin{gather*}
                         E_0\times M\xrightarrow{\sim} E\\
                         (\xi, m)\mapsto\Phi (m_0, m)\xi.
                \end{gather*}
                More invariantly, one identifies $ E_0$ with the
                space of $\nabla$-parallel sections of $ E$.
        \end{conditions}
\end{remarks}

\subsection{Dual connections and torsion}\lab{dual}
Let $\mathfrak g$ be a Lie algebroid and $\nabla$ a ${\mathfrak
  g}$-connection on itself. Generalizing \ref{abs}, we define the
{\df dual} of $\nabla$ to be the ${\mathfrak g}$-connection $\nabla^*$
on ${\mathfrak g}$ defined by
\begin{equation*}
        \nabla^*_XY:=\nabla_YX+[X, Y]_{\mathfrak g}.
\end{equation*}
One has `duality' in the sense that $\nabla^{**}=\nabla$.

The {\df torsion} of $\nabla$ is the section, $\torsion\nabla$, of
$\alternating^2 ({\mathfrak g},{\mathfrak g})$ measuring the
difference between $\nabla$ and its dual:
\begin{equation*}
        \torsion\nabla\, (X, Y):=\nabla_XY-\nabla^*_XY
        =\nabla_X Y-\nabla_YX-[X, Y].
\end{equation*}

The torsion or curvature of $\nabla$ can be expressed in terms  of 
the torsion and curvature of $\nabla^*$ (and, by duality, vice versa):
\begin{proposition}
        Let $\nabla$ be a ${\mathfrak g}$-connection on ${\mathfrak g}$, and 
        $\nabla^*$ its dual. Then:
        \begin{gather}
           \torsion\nabla=-\torsion\nabla^*\\
           \left\{
           \begin{split}
             \curvature\nabla\,(X, Y)Z&=(\nabla_Z^*\torsion\nabla^*)(X, Y)
                       +\curvature\nabla^*\,(X, Z)Y\\
                       &+\curvature\nabla^*\,(Z, Y)X;
                       \qquad X,Y,Z\subset\mathfrak g.
           \end{split}\right.\mathlab{dual2}
        \end{gather}
        In particular:
        \begin{conditions}
        \item\lab{scorch} If $\nabla^*$ is flat, then $\nabla$ is flat
          if and only if $\torsion\nabla^*$ is $\nabla^*$-parallel.
        \end{conditions}
\end{proposition}

\section{Cartan connections}\label{section5}%
In \ref{geom} a Cartan connection on a Lie algebroid ${\mathfrak
  g}$ was defined as a $TM$-connection on ${\mathfrak g}$ `compatible'
with the bracket on ${\mathfrak g}$.  Pivotal to understanding the
compatibility condition is the observation that the first jet bundle
$J^1{\mathfrak g}$ is also a Lie algebroid, a fact noted in
\cite{Higgins_Mackenzie_90}.  While the bracket on a Lie {\em algebra}
${\mathfrak g}$ determines a canonical representation
$\adjoint\colon{\mathfrak g}\rightarrow\gl ({\mathfrak g})$, the {\df
  adjoint} representation, the analogue for a Lie algebroid
${\mathfrak g}$ is a representation $\adjoint\colon J^1{\mathfrak
  g}\rightarrow\gl({\mathfrak g})$ defined by
$\adjoint_{J^1X}Y:=[X,Y]_{\mathfrak g}$. We show that a Cartan
connection $\nabla$ on ${\mathfrak g}$ is compatible with the bracket
on ${\mathfrak g}$ when the associated splitting $ J^1{\mathfrak
  g}\xleftarrow{s}{\mathfrak g}$ of
\begin{equation*}
        0\rightarrow T^*M\otimes{\mathfrak g}\monomorphism J^1{\mathfrak g}
        \rightarrow{\mathfrak g}\rightarrow 0
\end{equation*}
is a Lie algebroid morphism.  In that case, the
composite
\begin{equation*}
 {\mathfrak g}\xrightarrow{s} J^1{\mathfrak g}
 \xrightarrow{\adjoint}\gl({\mathfrak 
g})  
\end{equation*}
delivers a representation of $\mathfrak g$ on itself. There is
moreover, a representation of $\mathfrak g$ on $TM$, with respect to
which the anchor $\#\colon\mathfrak g\rightarrow TM$ is equivariant. Details
now follow.

\subsection{The first jet of a Lie algebroid}\lab{structure}
If ${\mathfrak g}$ is the Lie algebroid of a Lie groupoid $G$, then $
J^1{\mathfrak g}$ is a model for the Lie algebroid of $ J^1 G$, the
groupoid of one-jets of bisections of $G$.  This suggests a definition
for the bracket on $ J^1{\mathfrak g}$ more generally and this
definition agrees with the one to be described below. See also
\cite{Crainic_Fernandes_04}.

\begin{lemma}
        For any Lie algebroid $\mathfrak g$ over $M$, the formula
        \begin{equation*}
          \kappa_X(\alpha\otimes Y):=\mathcal L_{\#X}\alpha\otimes Y
          +\alpha\otimes[X,Y]_{\mathfrak g}
        \end{equation*}
        defines a Lie {\em algebra} representation,
        \begin{equation*}
          X\mapsto\kappa_X\colon \Gamma ({\mathfrak g})
          \rightarrow\gl\Big(\,\Gamma(T^*M\otimes{\mathfrak g})\,\Big).
        \end{equation*}
        Viewing sections of $ T^*M\otimes{\mathfrak g}$ as sections of
        $\homomorphism (TM,{\mathfrak g})$, one has
        \begin{equation}
                 (\kappa_X\phi)V:=[X,\phi (V)]_{\mathfrak g}+\phi 
                 ([V,\#X]_{TM}).\mathlab{alexia}
        \end{equation}
\end{lemma}

In making $ J^1{\mathfrak g}$ a Lie algebroid, the inevitable choice 
of anchor is the composite $ J^1{\mathfrak g}\rightarrow{\mathfrak 
g}\xrightarrow{\#}TM$, which we again denote by $\#$.  A Lie algebroid 
bracket on $ J^1{\mathfrak g}$ is a Lie bracket on the space of 
sections $ \Gamma (J^1{\mathfrak g})$, which, by Remark 
\eqrefs{TM}{TM1}, is isomorphic to $ \Gamma ({\mathfrak g})\oplus 
\Gamma (T^*M\otimes{\mathfrak g})$.  A natural choice of bracket is 
one making $ \Gamma (J^1{\mathfrak g})$ a semidirect product of 
$\Gamma({\mathfrak g})$ and $\Gamma (T^*M\otimes{\mathfrak g})$.

First, we make
$T^*M\otimes{\mathfrak g}$ a bundle of Lie algebras by
defining 
\begin{equation}
  [\phi_1,\phi_2]_{T^*M\otimes{\mathfrak
    g}}:=\phi_2\circ\#\circ\phi_1-\phi_1\circ\#\circ\phi_2.  \mathlab{bob}
\end{equation}
Next we use the representation $\kappa$ above to make $\Gamma
(J^1{\mathfrak g})$ a semidirect product of the Lie algebras
$\Gamma({\mathfrak g})$ and $\Gamma (T^*M\otimes{\mathfrak g})$. That
is, we define
\begin{align}
                [ J^1 X_1+\phi_1, J^1 X_2+\phi_2]_{J^1{\mathfrak g}}&:=J^1[ 
                X_1,X_2]_{\mathfrak g}+[\phi_1,\phi_2]_{T^*M\otimes{\mathfrak 
                g}}\mathlab{alex}\\
                &+\kappa_{X_1}\phi_2-\kappa_{X_2}\phi_1.\notag
\end{align}
In formulas \eqref{bob} and \eqref{alex}, $\phi_1,\phi_2\in \Gamma
(T^*M\otimes{\mathfrak g})$ and $ X_1, X_2\in \Gamma ({\mathfrak g})$
are arbitrary sections.  The compositions in \eqref{bob} make sense
because we identify sections of $T^*M\otimes\mathfrak g$ with
sections of $\homomorphism(TM,\mathfrak g)$.
        
With the above choice of anchor and bracket, $J^1\mathfrak g$ becomes
a Lie algebroid, containing $T^*M\otimes\mathfrak g$ as a subalgebroid. 

Our definition of the Lie algebroid structure of $ J^1 \mathfrak g$ is
natural in a categorical sense: Suppose that $ f\colon E\rightarrow F$
is a morphism of vector bundles over $ M$, covering the identity. Then
one defines $ J^1 f\colon J^1E\rightarrow J^1 F$ by $ J^1
f([\sigma]_m^1):=J^1 ([f\sigma]_m^1)$ and obtains the commutative
diagram
\begin{equation}
  \begin{CD}
        T^*M\otimes E @>>> J^1 E   @>>> E\\
        @V\alpha\otimes e\mapsto\alpha\otimes feVV   @VVJ^1fV  @VVfV\\
        T^*M\otimes F @>>>   J^1 F  @>>> F
  \end{CD}\qquad.\mathlab{diag}
\end{equation}
\begin{proposition}
         If $ f\colon{\mathfrak g}\rightarrow{\mathfrak h}$ is a morphism of 
         Lie algebroids, then so is $J^1f\colon J^1{\mathfrak g}\rightarrow 
         J^1{\mathfrak h}$.
\end{proposition}

\subsection{Adjoint representations}\lab{adjoint}
The Lie algebroid structure of $J^1\mathfrak g$ fixed above guarantees
the following result:
\begin{proposition}
For any Lie algebroid ${\mathfrak g}$, the formula
\begin{equation*}
        \adjoint_{J^1X}Y:=[X, Y]_{\mathfrak g},\qquad X,Y\in\Gamma(\mathfrak g),
\end{equation*}
defines a representation $\adjoint$ of $J^1\mathfrak g$ on $\mathfrak 
g$ that we call the {\df adjoint representation}.  If $\adjoint'\colon 
T^*M\otimes\mathfrak g\rightarrow\homomorphism(\mathfrak g,\mathfrak 
g)$ is defined by 
\begin{equation}
  \adjoint'(\alpha\otimes x)y:=-\alpha(\#y)x,\mathlab{giddy}
\end{equation}
then the following is a commutative diagram of Lie algebroid morphisms
with exact rows:
  \begin{equation*}
  \begin{CD}
        0 @>>> T^*M\otimes{\mathfrak g} @>\text{inclusion}>> 
               J^1{\mathfrak g}         @>>>                 {\mathfrak g}
      @>>> 0\\
      @VVV    
      @V\adjoint'VV
          @VV\adjoint V                      @VV\#V @VVV\\
        0 @>>> \homomorphism ({\mathfrak g},{\mathfrak g})
                                        @>\text{inclusion}>>
                                   \gl({\mathfrak g})       @>\#>> TM @>>> 0
  \end{CD}
  \end{equation*}
  For any smooth manifold $M$, the adjoint representation
  $\adjoint\colon J^1(TM)\rightarrow\gl(TM)$ is a Lie algebroid {\em
    isomorphism.}
\end{proposition}
\begin{caution}
  Note the minus sign in \eqref{giddy}, necessary for diagram
  commutativity. Note also the differing order of terms in definition
  \eqrefs{structure}{bob} of $[\phi_1,\phi_2]_{T^*M\otimes{\mathfrak
      g}}$, compared with that for $[\phi_1,\phi_2]_{\homomorphism
    ({\mathfrak g},{\mathfrak g})}$ in \eqrefs{glE}{cand} (taking
  $E=\mathfrak g$).\footnote{For this reason, in the special case
    ${\mathfrak g}=TM$, one must understand the canonical
    identification $T^*M\otimes TM\cong\homomorphism (TM,TM)$ as an
    {\em anti}homomorphism of Lie algebroids. This sign annoyance can
    only be shifted, not eliminated.}
\end{caution}

\subsection{Bracket compatibility}\lab{pat}
Let $\nabla$ denote an arbitrary $TM$-connection on a Lie algebroid
${\mathfrak g}$, let $ J^1{\mathfrak g}\xleftarrow{s}{\mathfrak g}$
denote the corresponding splitting of $0\rightarrow T^*M\otimes{\mathfrak
  g}\monomorphism J^1{\mathfrak g}\rightarrow{\mathfrak g}\rightarrow 0$
(described in \ref{TM}), and consider the composite morphisms
\begin{align}
        {\mathfrak g}\xrightarrow{s} J^1{\mathfrak g}\xrightarrow{\adjoint}\gl({\mathfrak 
        g})&,\mathlab{z1}\\
        {\mathfrak g}\xrightarrow{s} J^1{\mathfrak g}\xrightarrow{ 
        J^1\#} J^1(TM)\xrightarrow{\adjoint}\gl(TM)&.\mathlab{z2}
\end{align}
These define, respectively, a ${\mathfrak g}$-connection on 
${\mathfrak g}$ and a ${\mathfrak g}$-connection  on $TM$.  Both  
connections will be denoted $\bar\nabla$. Using \eqrefs{TM}{cost} and 
diagram \eqrefs{structure}{diag}, one computes for $ X\in \Gamma 
({\mathfrak g})$,
\begin{align}
        \bar\nabla_XY=\nabla_{\#Y}X+[X,Y]_{\mathfrak g},\qquad Y\in 
        \Gamma ({\mathfrak g}),\mathlab{z3}\\
        \bar\nabla_XV=\#\nabla_VX+[\#X, V]_{TM},\qquad V\in 
        \Gamma (TM).\mathlab{z4}
\end{align}
With respect to these connections, the anchor $\#\colon{\mathfrak 
g}\rightarrow TM$ is connection-preserving:
\begin{equation}
        \#\bar\nabla_XY=\bar\nabla_X\#Y,\qquad X, Y\in \Gamma({\mathfrak g}).
        \mathlab{z5}
\end{equation}
\begin{proposition}
        The connection $\nabla$ is compatible with the bracket on ${\mathfrak 
        g}$ in the sense of Definition \ref{geom} (i.e., is a Cartan 
        connection) if and only if  the corresponding splitting $ 
        J^1{\mathfrak g}\xleftarrow{s}{\mathfrak g}$ above is a Lie 
        algebroid morphism. In that case, the ${\mathfrak g}$-connections on 
        ${\mathfrak g}$ and $ TM$ given  by \eqref{z3} and \eqref{z4} are 
        representations.
\end{proposition}
\begin{proof}
With the assistance of \eqrefs{TM}{cost} and 
\eqrefs{structure}{alex}, we compute
\begin{align*}
        \curvature s\, (X, Y)=&[sX, sY]_{J^1{\mathfrak g}}-s[X, Y]_{\mathfrak 
        g}\\
        =&[J^1 X, J^1 Y]_{J^1{\mathfrak g}}-[J^1 X,\nabladot Y]_{J^1{\mathfrak g}}
        -[\nabladot X, J^1 Y]_{J^1{\mathfrak g}}\\
        &+[\nabladot X,\nabladot Y]_{J^1{\mathfrak g}}-J^1[ X, Y]_{\mathfrak g}+
        \nabladot[ X, Y]_{\mathfrak g}\\
        =&-\kappa_X(\nabladot Y)+\kappa_Y(\nabladot X)\\
        &+(\nabladot Y)\circ\#\circ (\nabladot X)-
         (\nabladot X)\circ\#\circ (\nabladot Y)\\
        &+\nabladot[ X, Y]_{\mathfrak g}.
\end{align*}
Applying \eqrefs{structure}{alexia}, we  obtain
\begin{align*}
        \curvature s\, (X, Y)V=&
         -[X,\nabla_VY]_{\mathfrak g}-\nabla_{[V,\#X]_{TM}}Y
         +[Y,\nabla_VX]_{\mathfrak g}+\nabla_{[V,\#Y]_{TM}}X\\
         &+\nabla_{\#\nabla_VX}Y-\nabla_{\#\nabla_VY}X+\nabla_V
           [X, Y]_{\mathfrak g}\\
         =&\nabla_V[X, Y]_{\mathfrak g}-[\nabla_VX, Y]_{\mathfrak g}
          -[X,\nabla_VY]_{\mathfrak g}\\
          &-\nabla_{\bar\nabla_YV}X+\nabla_{\bar\nabla_XV}Y,
          \quad\text{by \eqref{z4}.}
\end{align*}
Thus $\curvature s=0$ if and only if \eqrefs{geom}{compat} holds,
which proves the first claim.  The second claim holds because the
morphisms in \eqref{z1} and \eqref{z2} are Lie algebroid morphisms
whenever $ s\colon{\mathfrak g}\rightarrow J^1{\mathfrak g}$ is a Lie
algebroid morphism.  (The other morphisms appearing in \eqref{z1} and
\eqref{z2} are Lie algebroid morphisms by Propositions \ref{structure} and
Proposition \ref{adjoint}.)
\end{proof}
 
\subsection{On the curvature of Cartan connections}\lab{elong}
Let $\nabla$ be a Cartan connection. Then
\begin{equation}
        \curvature\nabla\, (\#X,\#Y)Z=(\bar\nabla_Z\torsion\bar\nabla) (X, 
        Y),\qquad X, Y, Z\in \Gamma ({\mathfrak g}),\mathlab{abba}
\end{equation}
where $\bar\nabla$ denotes the representation of ${\mathfrak g}$ on 
itself determined by $\nabla$.
\begin{proof}[Proof of \eqref{abba}]
         From \eqrefs{pat}{z3} follows the identity 
         $\nabla_{\#X}Y=\bar\nabla_X^*Y$, where $\bar\nabla^*$ is the dual of 
         $\bar\nabla\colon{\mathfrak g}\rightarrow\gl({\mathfrak g})$  (see 
         \ref{dual}).  So
         \begin{equation*}
                \curvature\nabla\, (\#X,\#Y)Z=\curvature\bar\nabla^*\, (X, Y)Z=
                (\bar\nabla_Z\torsion\bar\nabla) (X, Y),
         \end{equation*}
         by Proposition \eqrefs{dual}{dual2} and the flatness of 
         $\bar\nabla$.
\end{proof}
Combining \eqref{abba} with Theorem A, we have:
\begin{corollary}
  A transitive Cartan algebroid $ ({\mathfrak g},\nabla)$ is
  locally symmetric if and only if $\torsion\bar\nabla$ is
  $\bar\nabla$-parallel, where $\bar\nabla$ denotes the corresponding
  representation of $\mathfrak g$ on itself.
\end{corollary}
\noindent The transitive case is discussed further in Sect.\,\ref{section6}.

\subsection{Invariant differential operators}\lab{invariant}%
Let $({\mathfrak g},\nabla)$ be a Cartan algebroid on $M$. Let a
{\df ${\mathfrak g}$-tensor} on $M$ be any section of some
representation $E$ of ${\mathfrak g}$. As $\nabla$ determines
representations of ${\mathfrak g}$ on $ TM$ and itself, all
differential forms and all differential ${\mathfrak g}$-forms
(sections of $\alternating^k({\mathfrak g})$, for some $k$) are
$\mathfrak g$-tensors.

There is a {\df fundamental operator} $D$ acting on $\mathfrak
g$-tensors: if $\tau$ is a section of $ E$ and $\rho\colon{\mathfrak
  g}\rightarrow\gl(E)$ is a representation, then $ D\tau\in \Gamma
({\mathfrak g}^*\otimes E)$ is defined by
\begin{equation*}
        \langle D\tau, X\rangle=\rho_X\tau,\qquad X\in \Gamma ({\mathfrak 
        g}).
\end{equation*}
Since $ D\tau$ is another ${\mathfrak g}$-tensor, the fundamental 
operator $ D$ can be iterated to obtain higher order differential 
operators. 

As is well-known, there is an {\df exterior derivative operator} $d$
acting on $ E$-valued ${\mathfrak g}$-forms, whenever $ (E,\rho)$ is a
representation of ${\mathfrak g}$:

\begin{gather*}
        d\colon \Gamma (\alternating^k({\mathfrak g}, E))\rightarrow
         \Gamma (\alternating^{k+1}({\mathfrak g}, E))\\
         d\theta (X_0,\ldots, X_k):=\sum_{i=0}^k(-1)^i
         \rho_{X_k}(\theta (X_0,\ldots,\hat X_i,\ldots,X_k))\\
         -\sum_{i<j}(-1)^{i+j+1}\theta ([X_i, X_j]_{\mathfrak g},
          X_0,\ldots,\hat X_i,\ldots,\hat X_j,\ldots, X_k).
\end{gather*}
One has $d^2=0$, following from the fact that $\curvature\rho=0$.

Let $\omega\in \Gamma (\alternating^1 ({\mathfrak g},{\mathfrak g}))$
denote the tautological form, $\omega (X)=X$. Then $D\omega=0$ and
$d\omega=\torsion\bar\nabla$, where $\bar\nabla$ is the representation
of ${\mathfrak g}$ on itself induced by $\nabla$. The exterior
derivative $ d$ can be expressed in terms of the fundamental operator
$ D$ according to
\begin{equation*}
         d\theta=\omega\wedge D\theta+\theta^{d\omega}.
\end{equation*}
Here we view $ D\theta$ as a ${\mathfrak g}^*\otimes E$-valued 
${\mathfrak g}$-form and the wedge implies a contraction ${\mathfrak 
g}\otimes ({\mathfrak g}^*\otimes E)\rightarrow E$. Also,
\begin{equation*}
        \theta^{d\omega}(X_0,\ldots, X_k):=\sum_{i<j}(-1)^{i+j+1}
        \theta (d\omega(X_i, X_j),
          X_0,\ldots,\hat X_i,\ldots,\hat X_j,\ldots, X_k),
\end{equation*}
unless $ k=0$, in which case $\theta^{d\omega}:=0$.

\section{Transitive Cartan algebroids}\lab{section6}
The remainder of the paper concerns {\em transitive} Cartan
algebroids.
\subsection{Cartan connections recharacterized}\lab{tr}
As we have observed above, every Cartan connection $\nabla$ on a
Lie algebroid $\mathfrak g$ determines a representation of $\mathfrak
g$ on itself. In the transitive case the Cartan connection is {\em
  uniquely} determined by this representation, as we now explain.

Let $\mathfrak g$ be a transitive Lie algebroid and let $\mathfrak
h\subset\mathfrak g$ denote the kernel of its anchor. So we have an
exact sequence
\begin{equation*}
  0\rightarrow\mathfrak h\monomorphism\mathfrak g
  \xrightarrow{\#}TM\rightarrow 0.
\end{equation*}
The bracket on $\mathfrak g$ determines a canonical representation
$\rho$ of $\mathfrak g$ on $\mathfrak h$:
\begin{equation*}
  \rho_XY=[X,Y]_{\mathfrak g},
  \qquad X\in\Gamma(\mathfrak g),\,Y\in\Gamma(\mathfrak h).  
\end{equation*}
If $\nabla$ is a Cartan connection on
$\mathfrak g$, the corresponding representation $\bar\nabla$ of
$\mathfrak g$ on itself is given by
\begin{equation}
  \bar\nabla_XY:=\nabla_{\#Y}X+[X,Y]_{\mathfrak g}.\mathlab{trtwo}
\end{equation}
This representation leaves the subbundle $\mathfrak h\subset\mathfrak g$
invariant, and its restriction to $\mathfrak h$ is just the canonical representation $\rho$ above.
According to \eqref{trtwo},
\begin{equation*}
  \nabla_VX=\bar\nabla_XY+[Y,X]_{\mathfrak g}
  \quad\text{whenever}\quad\text{\#Y=V}.
\end{equation*}
The following result is readily verified:
\begin{proposition}
  Let $\mathfrak g$ be a transitive Lie algebroid with kernel
  $\mathfrak h$. Let $\mathfrak g\xleftarrow{t}TM$ be an arbitrary
  vector bundle splitting of the exact sequence 
  $$0\rightarrow\mathfrak h\rightarrow\mathfrak
  g\xrightarrow{\#}TM\rightarrow 0.$$
  Then for each representation $D$
  of $\mathfrak g$ on itself extending the canonical representation
  $\rho$ of $\mathfrak g$ on $\mathfrak h$, the $TM$-connection
  $\nabla$ on $\mathfrak g$, defined by
  \begin{equation*}
    \nabla_VX=D_XtV+[tV,X]_{\mathfrak g},
  \end{equation*}
  is independent of the choice of splitting $t$, and is a Cartan
  connection satisfying $\bar\nabla=D$. In particular, the map
  $\nabla\mapsto\bar\nabla $ is a bijection from the set of Cartan
  connections on $\mathfrak g$ to the set of representations of
  $\mathfrak g$ on itself extending $\rho$.
\end{proposition}

Given an arbitrary splitting $\mathfrak g\xleftarrow{t}TM$, we can also write \eqrefs{elong}{abba} as
\begin{equation}
   \curvature\nabla\, (V,W)Z= (\bar\nabla_Z\torsion\bar\nabla) (tV,tW ),
   \qquad V,W\in\Gamma(TM).\mathlab{hodge}
\end{equation}

\subsection{Reductive connections}\lab{red}
A Cartan connection $\nabla$ on $\mathfrak g$ will be called {\df
  reductive} if the corresponding representation $\bar\nabla$ of
$\mathfrak g$ on itself splits over the exact sequence $\mathfrak
h\rightarrow\mathfrak g\xrightarrow{\#}TM$. That is, if there exists a
splitting $\mathfrak g\xleftarrow{t}TM$ such that the induced
isomorphism $\mathfrak g\cong TM\oplus\mathfrak h$ of vector bundles
is an isomorphism of $\mathfrak g$-representations; it suffices to
check that $t$ is equivariant, i.e., that
\begin{equation}
  \bar\nabla_XtV=t\bar\nabla_XV,\qquad 
  X\in\Gamma(\mathfrak  g),\,\mathlab{bass}
  V\in\Gamma(TM),
\end{equation}
where $\bar\nabla$ on the left refers to the representation of
$\mathfrak g$ on itself determined by $\nabla$, while on the right it
refers to that on $TM$  (see \ref{pat}). One can express
the Cartan connection in terms of $t$ and $\bar\nabla$:
\begin{equation}
  \nabla_VX=t\bar\nabla_XV+[tV,X]_{\mathfrak g}.\mathlab{top}
\end{equation}
Conversely, for any splitting $t$ and representation $\bar\nabla$ of
$\mathfrak g$ on $TM$, \eqref{top} defines a reductive Cartan
connection.

\subsection{Infinitesimal $G$-structures}
A {\df $G$-structure} on a manifold $M$ is a reduction of its general
linear frame bundle; see \cite{Kobayashi_72}.  The corresponding
groupoid notion is a (transitive) subgroupoid of the frame groupoid
$\operatorname{GL}(TM)$, frame groupoids having been defined in
\ref{glE}.  Noting that the adjoint representation $\adjoint\colon J^1
(TM)\rightarrow\gl (TM)$ is an isomorphism, we define an {\df
  infinitesimal $G$-structure} to be any Lie subalgebroid of $J^1
(TM)$.\footnote{These are a special case of the {\df infinitesimal
    geometric structures} defined in \cite{Blaom_A}.} If ${\mathfrak
  g}\subset J^1 (TM)$ is an infinitesimal $G$-structure, then the
restriction $\adjoint\colon{\mathfrak g}\rightarrow\gl(TM)$ of the
adjoint representation (see \ref{adjoint}) is a canonical
representation of ${\mathfrak g}$ on $ TM$.  It is natural to restrict
attention to Cartan connections $\nabla$ on ${\mathfrak g}$ for which
the corresponding representation $\bar\nabla$ of ${\mathfrak g}$ on $
TM$ is this restricted adjoint representation.

The following generalizes Corollary \ref{Riemannian}.
\begin{theoremB}\lab{here}
  Consider a transitive infinitesimal $G$-structure ${\mathfrak
    g}\subset J^1 (TM)$ and let $\mathfrak h\subset T^*M\otimes TM$
  denote the kernel of its anchor.  Let $\nabla$ be a reductive Cartan
  connection on ${\mathfrak g}$ with $\bar\nabla=\adjoint$, and let
  ${\mathfrak g}\xleftarrow{ t} TM$ denote a corresponding splitting
  of $0\rightarrow\mathfrak h\monomorphism\mathfrak g\rightarrow
  TM\rightarrow 0$ (i.e., \eqrefs{red}{top} holds).  Let $\levi$
  denote the composite
  \begin{equation*}
         TM\xrightarrow{ t}{\mathfrak g}\xrightarrow{\adjoint}\gl(TM),
  \end{equation*}
  which is a $ TM$-connection on $TM$. Then the Cartan algebroid
  $({\mathfrak g},\nabla)$ is locally symmetric if and only if
  $\torsion\levi$ and $\curvature\levi$ are both ${\mathfrak
    h}$-invariant and $\levi$-parallel.
\end{theoremB}
\begin{proof}
  The curvature of $\nabla$ is given by \eqrefs{tr}{hodge}. With the
  help of the equivariance property \eqrefs{red}{bass}, one easily
  rewrites this as
  \begin{equation}
        \curvature\nabla\, (V,W)Z=\bar\nabla_Z(t^*\torsion\bar\nabla) (V, 
        W),\mathlab{f1}
  \end{equation}
  where
  \begin{equation*}
         t^*\torsion\bar\nabla\, (V, W):=\torsion\bar\nabla (tV, tW).
  \end{equation*}
  Appealing to equivariance and the definition of torsion, one obtains
  \begin{align}
         t^*\torsion\bar\nabla\, (V, W)&=t\bar\nabla_{tV}W-t\bar\nabla_{tW}V
         -[tV, tW]_{\mathfrak g}\notag\\
         &=t(\adjoint_{tV}W-\adjoint_{tW}V-[V, W]_{TM})-\curvature t\, (V, 
         W)\notag\\
         &=t(\torsion\levi\, (V, W))-\curvature t\, (V, W).\mathlab{f2}
  \end{align}
  We claim that
  \begin{equation}
        \curvature t\, (V, W)=-\curvature\levi\, (V, W).\mathlab{f3}
  \end{equation}
  Indeed, we have
  \begin{align*}
        \curvature\levi\, (V, W)&=\curvature (\adjoint\circ t)(V, W)\\
        &=([\adjoint (tV),\adjoint (tW)]_{\gl(TM)}-\adjoint (t[V,W]_{TM}))\\
        &=\adjoint ([tV, tW]_{\mathfrak g}-t[V, W]_{TM})\\
        &=\adjoint' (\curvature t\,(V, W)),
  \end{align*}
  implying (see \eqrefs{adjoint}{giddy}) that
  \begin{equation*}
        \curvature\levi\,(V, W) U=-\curvature t\,(V, W) U.
  \end{equation*}
  
  Substituting \eqref{f3} into \eqref{f2}, we obtain
  \begin{equation*}
         t^*\torsion\bar\nabla\,(V, W)=
         t(\torsion\levi\,(V, W))+\curvature\levi\,(V, W).
  \end{equation*}
  Using this and equivariance in \eqref{f1}, we arrive at the 
  formula
  \begin{equation*}
        \curvature\nabla\,(V, W) Z=t((\bar\nabla_Z\torsion\levi)(V, W))
        + (\bar\nabla_Z\curvature\levi)(V, W).
  \end{equation*}
  The first term on the right-hand side is a section of $ 
  t(TM)\subset{\mathfrak g}$, while the second term is a section of 
  ${\mathfrak h}$. So $\curvature\nabla=0$ if and only if 
  $\torsion\levi$ and $\curvature\levi$ are both $\bar\nabla$-parallel. 
  Now for any vector field $ V$ on $M$, we have
  \begin{align*}
        \bar\nabla_ZV&=Z(V)\quad\text{for}\quad Z\in \Gamma ({\mathfrak 
        h})\subset \Gamma (T^*M\otimes  TM),\\
        \text{while}\quad\bar\nabla_{tW}V&=\levi_WV\quad\text{for}\quad W\in 
        \Gamma (TM).
  \end{align*}
  So being $\bar\nabla$-parallel is the same as being simultaneously 
  ${\mathfrak h}$-invariant and $\levi$-parallel. The conclusion of 
  the theorem is now a consequence of Theorem A.
\end{proof}

\section{Cartan geometries}\lab{section7}
  Let $G_0$ be a Lie group and $H_0$ a closed subgroup.  Let 
  ${\mathfrak g}_0$ and ${\mathfrak h}_0$ denote the Lie algebras of 
  these groups.  Then a {\df Cartan geometry} modeled on $G_0/H_0$ is 
  a smooth manifold $M$, together with a principal $ H_0$-bundle 
  $\pi\colon P\rightarrow M$, and an absolute parallelism 
  $\omega\in\Omega^1 (P,{\mathfrak g}_0)$ of $ P$.  One requires that 
  $\omega$ is $H_0$-equivariant, where $ H_0$ acts on ${\mathfrak 
  g}_0$ via the adjoint action.  We  suppose $H_0$ acts on $P$ 
  from the right.  For details, see \cite{Sharpe_97}.
  
  A one-form $\omega$ as above is known as a {\df Cartan connection}, 
  but to distinguish it from the connections defined in \ref{geom}, we 
  will call $\omega$ a {\df classical Cartan connection}.

\subsection{Curvature} \lab{u}
The {\df curvature} of the Cartan geometry is the $\mathfrak
g_0$-valued two-form on $P$ defined by
\begin{equation}
  \Omega:=d\omega+\frac{1}{2}[\omega,\omega]_{\mathfrak g_0}.\mathlab{fad}
\end{equation}
The model $G_0/H_0$ is itself a Cartan geometry (take $ M=G_0/H_0$,
$P=G_0$, and let $\omega$ be the left-invariant Mauer-Cartan form on
$G_0$) and in this case the curvature vanishes.  Conversely, if $M$ is
an arbitrary Cartan geometry modeled on $G_0/H_0$, then $\Omega=0$
implies that $M$ is isomorphic to $G_0/H_0$, assuming certain global
obstructions also vanish \cite[Theorem 5.3]{Sharpe_97}. So curvature
is the local measure of deviation from the prescribed model.

We now describe a canonical, transitive, Cartan algebroid
$({\mathfrak g},\nabla)$ associated with a Cartan geometry $(M,
P\xrightarrow{\pi}M,\omega)$ modelled on $G_0/H_0$.

\subsection{The Lie algebroid  ${\mathfrak g}$}\lab{upi}%
The Lie algebroid is given by ${\mathfrak g}:=(TP)/H_0$, its anchor
being the map $\#\colon{\mathfrak g}\rightarrow TM$ sending $ v
\modulo H_0$ to $T\pi\cdot v$.  Sections of ${\mathfrak g}$ are in
one-to-one correspondence with the $H_0$-invariant vector fields on $
P$. Since the collection of such vector fields is closed under the
Jacobi-Lie bracket $[\,\cdot\,,\,\cdot\,]_{TP}$, we obtain a Lie
bracket $[\,\cdot\,,\,\cdot\,]_{\mathfrak g}$ on $ \Gamma ({\mathfrak
  g})$.  This bracket makes ${\mathfrak g}$ into a Lie algebroid.

An arbitrary vector field $X$ on $P$ is $H_0$-invariant if and only if
$\omega(X)$ is $H_0$-equivariant.

\subsection{The Cartan connection $\nabla$ on ${\mathfrak g}$} 
As an absolute parallelism, the classical Cartan connection $\omega$ 
determines a flat affine connection $D$ on $TP$.  Implicitly, $D$ is 
defined by
\begin{equation}
  \omega (D_XY)={\mathcal L}_X(\omega (Y)),\qquad 
  X,Y\in\Gamma(P),\mathlab{d1}
\end{equation}
where ${\mathcal L}$ is Lie derivative. On the other hand, if $X$ and
$Y$ are both $H_0$-invariant, then so is $D_XY$, which shows that $D$
may also be viewed as a flat $\mathfrak g$-connection on $\mathfrak
g$, i.e., as a representation of $\mathfrak g$ on itself.  By
Proposition \ref{tr}, $D$ determines a unique Cartan connection
$\nabla$ on $\mathfrak g$ satisfying $\bar\nabla=D$.

The Cartan algebroid $(\mathfrak g,\nabla)$ above encodes all local
information concerning the Cartan geometry, with the exception of
certain model information. One can recover the group $H_0$, up to
cover, but at best reconstructs $\mathfrak g_0$ as a representation of
$H_0$; unless $\nabla$ is flat (see Theorem C below), the Lie algebra
structure of $\mathfrak g_0$ is lost. We suppress details.

\subsection{Local symmetry}
For convenience, we now identify the $\mathfrak g_0$-valued two-forms
$d\omega$, $[\omega,\omega]_{\mathfrak g_0}$ and $\Omega$ with
$TP$-valued two-forms on $P$. This we may do using the absolute
parallelism $\omega$. From the definition of torsion it follows that $
\torsion D=d\omega$. It is straightforward to check that
$[\omega,\omega]_{\mathfrak g_0}$ is $D$-parallel. Equation
\eqrefs{u}{fad} then gives
\begin{equation*}
  D_Z\torsion D=D_Z\Omega,
  \qquad Z\in\Gamma(TP).
\end{equation*}
Specializing Corollary \ref{elong} to this case, we obtain:
\begin{theoremC}\lab{thmC}
   Let $ \pi\colon P\rightarrow M$ be a Cartan geometry with arbitrary 
   model data, and $\omega$ the corresponding classical Cartan 
   connection.  Let $ D$ denote the corresponding flat affine 
   connection on $ P$.  Let $ ({\mathfrak g},\nabla)$ denote the 
   corresponding Cartan algebroid described above.  Then $ ({\mathfrak 
   g},\nabla)$ is locally symmetric (in the sense of 
   \ref{dissymmetry}) if and only if the curvature $\Omega$ of the 
   Cartan geometry, viewed as a $TP$-valued two-form on $P$, is 
   $D$-parallel.
\end{theoremC}
%
\bibliography{dynamics2}
\end{document}